\documentclass{article}
\usepackage{amsmath,amssymb}

\newtheorem{thtype}{TheoremType}  
\newtheorem{definition}[thtype]{Definition}
\newtheorem{proposition}[thtype]{Proposition}
\newtheorem{example}{Example}

\newtheorem{rem}[thtype]{Remark}

\def\Rset{\mathbb{R}}

\def\proof{\noindent\mbox{\bf{Proof.}\ }}
\def\qed{\ \hfill \mbox{$\Box$}}

\begin{document}

\title{Existence of an observation window of finite width for continuous-time autonomous nonlinear systems
    \thanks{
      This work was supported by the Japan Society for the Promotion of 
      Science under Grant-in-Aid for Scientific Research (C) 23560535. 
      This document is the accepted version of the manuscript published in Automatica 
      in Volume 75 with 
      DOI: https://doi.org/10.1016/j.automatica.2016.08.005
      and is a replacement of the preparatory version given in arXiv,
      which the author uploaded at the time instant of its submission.  
      If possible, please use the official version instead of this file.  
      $\copyright$2018. This manuscript version is made available under the CC-BY-NC-ND 4.0 license 
      http://creativecommons.org/licenses/by-nc-nd/4.0/       
    }
  }
\author{Shigeru Hanba
      \thanks{
      Department of Electrical and Electronics Engineering, University of the Ryukyus,
      1 Senbaru Nishihara, Nakagami-gun, Okinawa 903-0213, Japan; email: sh@gargoyle.eee.u-ryukyu.ac.jp
    }}
\maketitle

  \begin{abstract}
    In this note, the relationship between notions of observability
    for continuous-time nonlinear system related to distinguishability, 
    observability rank condition and K-function 
    has been investigated.  It is proved that an autonomous nonlinear system
    that is observable in both distinguishability and rank condition sense 
    permits an observation window of finite width, and 
    it is possible to construct a K-function related to observability for such system.  
  \end{abstract}

\section{Introduction}
The state estimation problem is one of the most fundamental problems in control system theory.  
The intrinsic property of the system that makes the state estimation problem feasible is
observability, and this property has been extensively studied in past decades.  

For linear systems, the notion of observability is firmly established\cite{Kailath1980}.  
Contrary, for nonlinear systems, several non-equivalent definitions of observability have been proposed,
and by using these definitions, the state estimation problem of nonlinear systems have been extensively studied
\cite{Alamir2007,Besancon2007,Gauthier1992,Hammouri2007,Hermann1977,Isidori1995,Nijmeijer1990}.  
However, the relations between several different notions of observability have not been fully understood
(although there are several established facts\cite{Besancon2007,Hammouri2007,Hermann1977}.)  

Recently, the author has proved that
a discrete-time nonlinear systems whose observation map 
is injective and the Jacobian of the observation map is of full rank
satisfies a seemingly stronger condition called uniform observability\cite{Hanba2009},
and by assuming these conditions, it is possible to construct a
${\mathcal K}$-function related to observability\cite{Hanba2010}.  
This note is an attempt to establish corresponding results for continuous-time systems.  

The scope of this note is limited to autonomous nonlinear systems,
and we deal with three typical definitions of observability for nonlinear systems
which are related to distinguishability, rank condition and ${\mathcal K}$-function, respectively
(precise definitions are given later.)  
Roughly speaking, we prove that distinguishability together with the observability rank condition
implies that there is a `observation window' (the sequence of past output as a function of time)
of finite width which determines the initial state uniquely,
and it is possible to construct a ${\mathcal K}$-function related to observability.  

A preliminary version of this manuscript is available in arXiv\cite{Hanba2015}.  

\section{Main Results}
In this note, we consider an autonomous nonlinear system of the form 
\begin{equation}
  \begin{split}
    \dot{x}&=f(x), \\
    y&=h(x), 
  \end{split}
  \label{eq:a_nls_01}
\end{equation}
where $x \in \Rset^{n}$ is the state and $y \in \Rset^{p}$ is the output.  
The functions $f(x)$ and $h(x)$ are assumed to be of compatible dimensions, 
and smooth up to required order.  The solution of  (\ref{eq:a_nls_01}) is assumed to be unique,
and the solution initialized at $t=0$ by $x_0$ is denoted by $\varphi(t,0,x_0)$, 
which is assumed to be a continuous function of $x_0$.  
The set of permissible initial conditions of (\ref{eq:a_nls_01}) is assumed to be compact,
and is denoted by $\Omega$.  Note that the state space itself is not necessarily compact.  

Next, we introduce the definitions of observability considered in this note.  
Unfortunately, there is no general agreement on how to name these properties.  
Hence, we temporally call them `D-observability', 'R-observability', and `K-observability'
(to be defined below) for brevity.  

\begin{definition}\label{def:indistinguishable}\cite{Besancon2007}
  A pair of initial states $(x_1,x_2)$ of (\ref{eq:a_nls_01}) 
  with $x_1 \neq x_2$ is said to be an indistinguishable pair if
  $\forall t \geq 0$, $h(\varphi(t,0,x_1)) = h(\varphi(t,0,x_2))$.  
\end{definition}

\begin{definition}\label{def:observable}\cite{Besancon2007}
  The system (\ref{eq:a_nls_01}) is said to be D-observable 
  (with respect to $\Omega$)
  if there is no indistinguishable pair in the set $\Omega$.  
\end{definition}

\begin{definition}\label{def:rank-observable}\cite{Besancon2007}
  The system (\ref{eq:a_nls_01}) is said to be R-observable (with respect to $\Omega$)
  if $\exists N>0$, the Jacobian of the map $H(x)=(h(x),L_f h(x),\ldots,L_f^{N-1}h(x))$ is
  of full rank on $\Omega$, 
  where $L_f h=\frac{\partial h}{\partial x}f$,
  and $L_f^{k} h=L_f \left ( L_f^{k-1} h \right )$.  
\end{definition}

\begin{definition}\cite{Haddad2008,Khalil1996}
  A function $\alpha: D \rightarrow [0,\infty)$ 
  (where $D$ is either $[0,\infty)$, $[0,a)$ or $[0,a]$ with $a>0$)
  is said to be  a ${\mathcal K}$-function if it is continuous, $\alpha (0)=0$, 
  and is strictly increasing. 
\end{definition}

\begin{definition}\label{def:K-observable}\cite{Alamir2007}
  The system (\ref{eq:a_nls_01}) is said to be K-observable
  (with respect to $\Omega$) if  
  $\exists T>0$, $\forall x_1, x_2 \in \Omega$, 
  \begin{equation}
    \int_{0}^{T} \left | h(\varphi(t,0,x_1))-h(\varphi(t,0,x_2)) \right |^2 dt \geq \alpha(|x_1-x_2|), 
    \label{eq:K_observable_01}
  \end{equation}
  where $\alpha(\cdot)$ is a ${\mathcal K}$-function and $|\cdot|$ denotes the Euclid norm of a vector.   
\end{definition}

\begin{rem}
Each of above definitions require smoothness of $f(x)$ and $h(x)$ in different level.
\begin{itemize}
\item For Definition~\ref{def:observable}, the only requirement is that the system (\ref{eq:a_nls_01}) 
  has a unique solution. A finite escape time is allowed, as far as the state distinction is 
  achievable before the arrival of the finite escape time.  There is no restriction to $h(x)$.  
\item For Definition~\ref{def:rank-observable}, $h(x)$ should be $N-1$ times continuously differentiable,
  and $f(x)$ should be $N-2$ times continuously differentiable, but the value of $N$ cannot be specified
  (although it is finite).  
\item For Definition~\ref{def:K-observable}, the requirements are that (\ref{eq:a_nls_01}) 
  has a unique solution, the solution of (\ref{eq:a_nls_01}) is defined for $t \in [0,T]$, 
  and $h(\varphi(t,0,x_0))$ is integrable for each $x_0 \in \Omega$.  
\end{itemize}
\end{rem}

It is a known fact that, if a system is R-observable at a point $x_0$, 
then it is D-observable  on a neighborhood of $x_0$\cite{Besancon2007,Hermann1977},
but it is not always possible to extend the result to the whole of $\Omega$.  
On the other hand, D-observability does not imply R-observability, as the following example shows.  
\begin{example}
  Consider a 1-dimensional system
  \begin{equation}
    \begin{split}
      \dot{x}&=-x,\\
      y&=h(x)=x^3
    \end{split}  
    \label{eq:c_example02}
  \end{equation}
  This system is D-observable because it is possible to directly calculate $x$ from $y$
  ($x=y^{1/3}$),   but is not R-observable at $x=0$,
  because $h(x)=x^3$, $L_f h(x)=-3 x^3$, and inductively, 
  $L_f^k h(x)=(-1)^k 3^k x^{3}$,
  and hence their derivatives vanish at $x=0$.  
\end{example}

It is desirable that the width of the `observation window'
(the time interval that the output of the system is stored in order to determine the initial state uniquely)
is finite.  In this sense, K-observability is convenient,
and has been widely adopted in works on moving horizon state estimation\cite{Alamir2007,Alessandri2008}.  
If the system (\ref{eq:a_nls_01}) is K-observable, then for $x_1,x_2 \in \Omega$ with $x_1 \neq x_2$, 
$\exists t : 0 \leq t \leq T$, $h(\varphi(t,0,x_1)) \neq h(\varphi(t,0,x_2))$, hence 
(\ref{eq:a_nls_01}) is D-observable.  
Then, a natural question arises: do systems that are D-observable
always permit an observation window of finite width?  
Unfortunately, the answer is negative, which is given in the following example.  
\begin{example}
  Consider a 1-dimensional system
  \begin{equation}
    \begin{split}
      \dot{x}&=x\\
      y&=h(x)=
      \begin{cases}
        0 & x< M\\
        x-M & x \geq M, 
      \end{cases}
    \end{split}  
    \label{eq:c_example01}
  \end{equation}
  where $M$ is a positive constant.  
  If the initial condition is zero, then the output is identically zero.  
  For an initial condition $x_0 >0$, 
  $x(t)=\exp[t]x_0$, and hence the output is identically zero
  for $t < \ln (M/x_0)$ and is $\exp[t] x_0-M$ for $t \geq \ln (M/x_0)$.  
  Hence, the zero initial condition and $x_0$ cannot be distinguished until $t=\ln (M/x_0)$,
  and hence as the initial condition gets smaller, 
  the required width of the observation window tends to infinity.  

  One may argue that the reason for making the width of the observation window infinite 
  is the non-differentiability of the output function, but this is not the case.  
  For example, by replacing the output function $h(x)$ of (\ref{eq:c_example01}) with 
  \[
  \begin{split}
    h(x)&=
    \begin{cases}
      0 & x \leq M\\
      \exp[-1/(x-M)] & x \geq M
    \end{cases}, 
  \end{split}
  \]
  a similar conclusion holds.  
\end{example}

Thus far, we have seen that there are gaps between 
D-observability, R-observability and K-observability, 
and a D-observable system does not always permit an observation window of finite width.  
In the following, we show that,
if (\ref{eq:a_nls_01}) is D-observable as well as R-observable, 
then there exists an observation window of finite width, 
and it is possible to construct a ${\mathcal K}$-function corresponding to Definition~\ref{def:K-observable},
and hence (\ref{eq:a_nls_01}) is K-observable.

\begin{proposition}\label{th:f_o}
  If (\ref{eq:a_nls_01}) is D-observable as well as R-observable 
  for the initial condition set $\Omega$, then there is a finite $T>0$
  such that $\forall x_1,x_2 \in \Omega$ with $x_1 \neq x_2$,
  $\exists t: 0 \leq t \leq T$, $h(\varphi(t,0,x_1)) \neq h(\varphi(t,0,x_2)) $.  
\end{proposition} 

\proof
We first prove that 
\begin{equation}
  \begin{split}
  &\forall x \in \Omega, \exists {\mathcal N}(x), \forall z_1,z_2 \in {\mathcal N}(x) {\text{ with }} z_1 \neq z_2, \\
  &\quad \forall T>0, \exists t: 0 \leq t \leq T, h(\varphi(t,0,z_1)) \neq h(\varphi(t,0,z_2))     
  \end{split}
  \label{eq:finite_o_w_01}
\end{equation}
by contradiction, where ${\mathcal N}(x)$ denotes an open neighborhood of $x$.  
Suppose that (\ref{eq:finite_o_w_01}) is false, that is, 
\begin{equation}
  \begin{split}
  &\exists x \in \Omega, \forall {\mathcal N}(x), \exists z_1,z_2 \in {\mathcal N}(x) {\text{ with }}  z_1 \neq z_2, \\
  &\quad \exists T>0, \forall t: 0 \leq t \leq T, h(\varphi(t,0,z_1)) = h(\varphi(t,0,z_2)).      
  \end{split}
  \label{eq:finite_o_w_02}
\end{equation}
Then, $h(\varphi(t,0,z_1)) - h(\varphi(t,0,z_2))$ is identically zero as a function of $t$.  
Hence, for all $k \geq 0$, $\frac{d^k}{dt^k}h(\varphi(t,0,z_1)) =\frac{d^k}{dt^k}h(\varphi(t,0,z_2))$,
hence $H(\varphi(t,0,z_1))=H(\varphi(t,0,z_2))$, and by letting $t=0$, $H(z_1)=H(z_2)$.  
On the other hand, because the Jacobian of $H$ is of full rank, 
there is a neighborhood of $x$ in which $H$ is injective.  By choosing such neighborhood ${\mathcal N}(x)$
(recall that ${\mathcal N}(x)$ is arbitrary), it follows that $H(z_1) \neq H(z_2)$ because $z_1 \neq z_2$,
hence a contradiction has been obtained.  
Therefore, (\ref{eq:finite_o_w_02}) is false and hence (\ref{eq:finite_o_w_01}) is true.  

Next, fix $x \in \Omega$, and let ${\mathcal N}_{\rm{loc}}(x)$ be an open neighborhood of $x$ in which 
\begin{equation}
  \begin{split}
  &\forall z_1,z_2 \in {\mathcal N}_{\rm{loc}}(x) {\text{ with }} z_1 \neq z_2, \\
  & \quad \forall T>0, \exists t: 0 \leq t \leq T, h(\varphi(t,0,z_1)) \neq h(\varphi(t,0,z_2))     
  \end{split}
  \label{eq:finite_o_w_03}
\end{equation}
holds.  Because $\Omega$ is compact and ${\mathcal N}_{\rm{loc}}(x)$ is open,
$\Omega \setminus  {\mathcal N}_{\rm{loc}}(x)$ is compact.  
For each $z \in \Omega \setminus  {\mathcal N}_{\rm{loc}}(x)$, by D-observability, 
$\exists t_z >0, h(\varphi(t_z,0,x)) \neq h(\varphi(t_z,0,z))$.  
Because $h$ and $\varphi$ are continuous, there is an open neighborhood $O_z$ of $x$ and an open neighborhood $G_z$ of z
for which $h(\varphi(t_z,0,O_z)) \cap h(\varphi(t_z,0,G_z))=\emptyset$.  
Let $W_z=(z,t_z,O_z,G_z)$ be the tuple satisfying this condition, 
and consider the set
\begin{equation}
  \{W_z: z \in \Omega \setminus {\mathcal N}_{\rm{loc}}(x)\}.  
  \label{eq:f_o_tuple_01}
\end{equation}
Because $\{G_z: z \in \Omega \setminus {\mathcal N}_{\rm{loc}}(x)\}$ covers $\Omega \setminus {\mathcal N}_{\rm{loc}}(x)$
and $\Omega \setminus {\mathcal N}_{\rm{loc}}(x)$ is compact, 
for some $L>0$, there is a finite subcollection
\[
\{W_z^{(1)},\ldots,W_z^{(L)}\}
\]
of (\ref{eq:f_o_tuple_01}) (let us rewrite $W_z^{(i)}=(z^{(i)},t_z^{(i)},O_z^{(i)},G_z^{(i)})$)
for which $\{G_z^{(1)},\ldots,G_z^{(L)}\}$ covers $\Omega \setminus {\mathcal N}_{\rm{loc}}(x)$.  
Let 
$V_x =\left ( \cap_{i=1}^{L} O_z^{(i)}\right ) \cap {\mathcal N}_{\rm{loc}}(x)$, 
\begin{equation}
  T_x =\max\{ t_z^{(1)},\ldots,t_z^{(L)}\}.  
  \label{eq:f_o_max_time_01}
\end{equation}
$V_x$ is an open neighborhood of $x$, because it is a finite intersection of open sets containing $x$.  
Let $M_x=(x,V_x,T_x)$ be the tuple corresponding to above construction, 
and consider the set 
\begin{equation}
  \{M_x: x \in \Omega \}.  
  \label{eq:f_o_tuple_02}
\end{equation}
Because $\Omega$ is compact, 
for some $J>0$, 
there is a finite subcollection 
$\{M_x^{(1)},\ldots,M_x^{(J)}\}$ of (\ref{eq:f_o_tuple_02}) 
(let us rewrite $M_x^{(j)}=(x^{(j)},V_x^{(j)},T_x^{(j)})$)
for which $\{V_x^{(1)},\ldots,V_x^{(J)}\}$ covers $\Omega$.  
Let 
\begin{equation}
  T=\max\{ T_x^{(1)},\ldots,T_x^{(J)}\}.  
  \label{eq:f_o_max_time_02}
\end{equation}
Then, all pairs of initial conditions $(x_1,x_2)$ with $x_1 \neq x_2$ are distinguishable for some 
$t$ with $0 \leq t \leq T$.  
For, by construction, there is a $V_x^{(j)}$ such that $x_1 \in V_x^{(j)}$.  
If $x_2$ is in ${\mathcal N}_{\rm{loc}}(x^{(j)})$, then by (\ref{eq:finite_o_w_03}),
$x_1$ and $x_2$ are distinguishable.  Otherwise, 
$z_2 \in \Omega \setminus {\mathcal N}_{\rm{loc}}(x^{(j)})$. 
Because $\Omega \setminus {\mathcal N}_{\rm{loc}}(x^{(j)})$ is covered by corresponding $\{G_z^{(1)},\ldots,G_z^{(L)}\}$,
by (\ref{eq:f_o_max_time_01}) and (\ref{eq:f_o_max_time_02}),
$x_1$ and $x_2$ are distinguishable at some $t$ with $0 \leq t \leq T$.    
\qed

\begin{rem}
In Proposition~\ref{th:f_o}, it has been implicitly assumed that 
(\ref{eq:a_nls_01}) has a well-defined solution for $t \in [0,T]$.  
\end{rem}

Next, we prove that a system that is D-observable as well as R-observable on $\Omega$ is K-observable.  

\begin{proposition}\label{th:K}
  If the system (\ref{eq:a_nls_01}) is D-observable as well as R-observable on $\Omega$,
  then it is K-observable on $\Omega$.  
\end{proposition}

\proof
Instead of constructing a ${\mathcal K}$-function that satisfies (\ref{eq:K_observable_01}) directly,
we construct an increasing function of $|x_1-x_2|$ that is positive if $|x_1-x_2| \neq 0$.  
We assume that $T$ of (\ref{eq:K_observable_01}) is sufficiently large and
$\forall x_1,x_2 \in \Omega$ with $x_1 \neq x_2$, $\exists t: 0 \leq t \leq T$, 
$h(\varphi(t,0,x_1))\neq h(\varphi(t,0,x_2))$.  The existence of such $T$ is assured by Proposition~\ref{th:f_o}.  
By using this $T$, let
\[
\theta(x_1,x_2)=\int_{0}^{T}|h(\varphi(t,0,x_1))-h(\varphi(t,0,x_2))|^2 dt.  
\]
Note that $\theta(x_1,x_2)$ is a continuous function of $(x_1,x_2)$.  
By the construction of $T$, if $x_1 \neq x_2$, 
then $\exists t: 0 \leq t \leq T$ with $h(\varphi(t,0,x_1)) \neq h(\varphi(t,0,x_2))$, 
and hence $\theta(x_1,x_2) > 0 $.  
For $r>0$, let  $D_r = \{(x_1,x_2) \in \Omega \times \Omega: |x_1-x_2| \geq r\}$,
which is a compact set.  
Let $\alpha_0(r)=\min_{(x_1,x_2) \in \Xi(r)} \theta(x_1,x_2)$. Because $D_r$ is compact,
the minimum is well defined and 
because $(x_1,x_2) \in \Xi(r)$ implies that $x_1 \neq x_2$, the minimum is positive.  
If $r_1<r_2$, then $D_{r_1} \supset D_{r_2}$, hence $\alpha_0(r_1) \leq \alpha_0(r_2)$.  
Therefore, $\alpha_0(r)$ is a monotone nondecreasing function of $r$ and its value is positive.  
For definiteness, let $\alpha_0(0)=0$.  
Note that $\alpha_0(r)$ may not be strictly increasing, and may be discontinuous.  
However, by Lemma~5 of \cite{Hanba2010}, it is possible to construct a ${\mathcal K}$-function 
that satisfy $\alpha(r) \leq \alpha_0(r)$.  Hence, 
\begin{multline*}
\int_{0}^{T}|h(\varphi(t,0,x_1))-h(\varphi(t,0,x_2))|^2 dt\\
 \geq \alpha_0(|x_1-x_2|) \geq \alpha(|x_1-x_2|),     
\end{multline*}
whence (\ref{eq:a_nls_01}) is K-observable.  
\qed

\begin{rem}
  If all of $f(x)$, $h(x)$ and $\varphi(t,0,x_0)$ are real analytic functions of their arguments and the system is forward complete, 
  it can be proved that D-observability implies K-observability without assuming R-observability, 
  and the width of the observation window of Proposition~\ref{th:K} may be arbitrarily small.  
  The proof is as follows.  
  First, if $h(\varphi(t,0,x_0))$ is real analytic and $h(\varphi(t,0,x_1)) \equiv h(\varphi(t,0,x_2))$
  on an interval $[0,T)$ for a pair $(x_1,x_2)$ with $x_1 \neq x_2$, 
  then $h(\varphi(t,0,x_1)) \equiv h(\varphi(t,0,x_2))$ for any $t$ due to real analyticity.  
  Therefore, if (\ref{eq:a_nls_01}) is D-observable, then $\forall T>0$, 
  $\forall x_1,x_2$ with $x_1 \neq x_2$, $\exists t: 0 \leq t \leq T$, 
  $h(\varphi(t,0,x_1)) \neq h(\varphi(t,0,x_2))$.  
  Hence, for real analytic systems,
  the existence of the finite observation window of Proposition~\ref{th:f_o} is
  established without assuming R-observability, and its width may be arbitrarily small.  
  For any $T>0$, the construction of a K-function is identical to the proof of Proposition~\ref{th:K}.  
  Hence, we have proved that, if (\ref{eq:a_nls_01}) is D-observable and forward complete, 
  and all of $f(x)$, $h(x)$ and $\varphi(t,0,x_0)$ are real analytic functions of their arguments, 
  then (\ref{eq:a_nls_01}) is K-observable for any observation window of positive length.  
  Because a K-observable system is automatically D-observable (without assuming real analyticity), 
  we have shown that, for a forward complete and real analytic system,
  the notions of D-observability and K-observability are equivalent.  
\end{rem}

\section{Conclusion}
In this note, we have shown that
an autonomous nonlinear system which is D-observable and R-observable
always permits an observation window of finite width,
and it is actually K-observable as well.  
A theoretical construction of corresponding ${\mathcal K}$-function has been provided as well.  
In many researches related to nonlinear observability, 
the existence of an observation window of finite width and a ${\mathcal K}$-function
are assumed {\em{a priori}}.  Contrary, in this note, it has been proved that 
they are consequences of D-observability and R-observability.  
It is also to be noted that our result is purely existential, 
and no practical method of obtaining the width of the observation window 
and the ${\mathcal K}$-function have been provided.


\begin{thebibliography}{99}

\bibitem{Alamir2007}
  M.~Alamir, 
  Nonlinear moving horizon observers: theory and real-time implementation.  
  In Besan\c{c}on G. (Ed.), 
  {\em{Nonlinear Observers and Applications}}, 
  Springer-Verlag, pp. 139--179, 
  2007.  

\bibitem{Alessandri2008}
  A.~Alessandri, M.~Baglietto and G.~Battistelli,   
  Moving-horizon state estimation for nonlinear discrete-time systems: new stability results and approximation scheme.  
  Automatica, 
  Vol.~44, No.~7, pp.~1753--1765,
  2008.  

\bibitem{Besancon2007}
  G.~Besan\c{c}on, 
  An overview on observer tools for nonlinear systems, 
  In G.~Besan\c{c}on (Ed.), 
  {\em{Nonlinear Observers and Applications}},  
  Springer-Verlag,  pp.~1--33, 
  2007.  

\bibitem{Gauthier1992}
  J.~P.~Gauthier, H.~Hammouri and S.~Othman,   
  A simple observer for nonlinear systems applications to bioreactors, 
  IEEE Transactions on Automatic Control, 
  Vol.~37, No.~6,  pp.~875--880, 
  1992.  

\bibitem{Haddad2008}
  W.~M.~Haddad and  V.~Chellaboina,   
  {\em{Nonlinear Dynamical Systems and Control: A Lyapunov-Based Approach}}, 
  Princeton University Press,
  2008.  

\bibitem{Hammouri2007}
  H.~Hammouri, 
  Uniform observability and observer synthesis, 
  In G.Besan\c{c}on (Ed.), 
  {\em{Nonlinear Observers and Applications}}, 
  Springer-Verlag, pp.~35--70,
  2007.  

\bibitem{Hanba2009}
  S.~Hanba,   
  On the `uniform' observability of discrete-time nonlinear systems, 
  IEEE Transactions on Automatic Control, 
  Vol.~54, No.~8, pp.~1925--1928, 
  2009.  

\bibitem{Hanba2010}
  S.~Hanba, 
  Further results on the uniform observability of discrete-time nonlinear systems, 
  IEEE Transactions on Automatic Control, 
  Vol.~55, No.~4, pp.~1034--1038, 
  2010.  

\bibitem{Hermann1977}
  R.~Hermann, and A.~Krener,  
  Nonlinear controllability and observability, 
  IEEE Transactions on Automatic Control, 
  Vol.~22, No.~5, pp.~728--740, 
  1977.  

\bibitem{Isidori1995}
  A.~Isidori,   
  {\em{Nonlinear Control Systems}},  3/e, 
  Springer-Verlag, London, UK, 
  1995.  

\bibitem{Khalil1996}
  H.~K.~Khalil,  
  {\em{Nonlinear Systems}}, 2/e, 
  Prentice-Hall, 
  1996.  

\bibitem{Kailath1980}
  T.~Kailath,  
  {\em{Linear Systems}}, 
  Prentice-Hall,
  1980.  

\bibitem{Nijmeijer1990}
  H.~Nijmeijer and A. J. van der Schaft,  
  {\em{Nonlinear Dynamical Control Systems}}, 
  Springer-Verlag, New York, USA, 
  1990.  

\bibitem{Hanba2015}
  S.~Hanba,  
  Existence of an observation window of finite width for continuous-time autonomous nonlinear systems, 
  arXiv:1505.05970,   2015.  
  http://arxiv.org/pdf/1505.05970v1

\end{thebibliography}
\end{document}